# Periodicity Properties of Coefficients of Half Integral Weight Modular Forms [*]

Alexandru Tupan
*McGill University/Concordia University/CICMA*

December 4, 2001

**Abstract.** Let $\sum_{n \geq 1} a(n) q^n$ be the $q$-expansion at $\infty$ of a half integral weight modular form $f$. We consider blocks of coefficients of $f$ of the form $\{a(dn^2) | n \in \mathbb{N}\}$ where $d$ is a square free natural integer. Under some weak assumptions on $f$, we show that if such a coefficient block is periodic when considered as a function of $n$, then it must vanish completely. The proof is analytic in nature and uses Shimura's lifting theorem together with estimates on the order of growth of Fourier coefficients of modular forms.

## 1. Introduction

Consider the space $S_{k+\frac{1}{2}}(N, \chi)$ of cusp forms of half-integral weight, level $N$, and character $\chi$. When the weight is $1/2$ a well known result (J.-P. Serre and H.M. Stark, 1976) states that any element of $S_{\frac{1}{2}}(N, \chi)$ is a linear combination of theta series of the form

$$\theta_{\psi,d} = \sum_{n \geq 1} \psi(n) q^{dn^2}.$$

In particular if $f = \sum_{n \geq 1} a(n) q^n$ is the Fourier expansion at the cusp at $\infty$ of a cusp form $f$ of weight $1/2$, then any block of coefficients of the form

$$a(d), a(4d), a(9d), \ldots, a(n^2 d), \ldots$$

is a periodic function in $n$.

A natural question to ask is whether a similar periodicity result remains true in more generality. In this paper we shall answer such a question in the case of cusp forms of weight at least $3/2$. We prove that, under very weak technical conditions on the form $f$, a periodic block of coefficients must be completely zero. The precise statement is given in section 3, Theorem A.

In the proof we shall essentially use two results of a relatively different nature: Shimura's Main Theorem (G. Shimura, 1973), and a result of Hecke regarding the modularity of a product of two Dirichlet $L$-series ((T. Miyake, 1989), §4.7). They will give us the information we need, when combined with estimates on the order of growth of coefficients of cusp forms (the integral weight case).

---

[*] E-mail: tupan@math.mcgill.ca

In a remark after the proof of Theorem A we shall point out how our method can be used to prove results similar in nature. As example we can show the vanishing of a block of coefficients for which the sequence $\{\frac{a(n^2 d)}{n} | n \in \mathbb{N}\}$ is periodic in $n$ (the periodicity of such a sequence is always true for $\theta$-series of weight $3/2$). The statement is true under the hypothesis that the weight is at least $5/2$.

## 2. Preliminaries

### 2.1. Variables.

We use standard notations. The letter $\mathbf{H}$ denotes the upper half-plane $\{z | \text{Im}(z) > 0\}$. $z$ will always be a complex number in the upper half plane. We put $q = e^{2\pi i z}$. The letter $s$ will denote a complex variable, $s \in \mathbb{C}$.

### 2.2. Spaces of forms.

Let $k$ be a positive integer and $\chi : (\mathbb{Z}/N\mathbb{Z})^\times \to \mathbb{C}^\times$ be a character (mod $N$). By $M_k(N, \chi)$ we denote the space of modular forms of weight $k$, level $N$, and Nebentypus character $\chi$. Similarly, when $4|N$, $M_{k+\frac{1}{2}}(N, \chi))$ denotes the space of modular forms of weight $k + \frac{1}{2}$, level $N$ and Nebentypus character $\chi$. We denote by $S_k(N, \chi)$ and $S_{k+\frac{1}{2}}(N, \chi)$ the corresponding subspaces of cusp forms.
$S^*_{\frac{3}{2}}(N, \chi)$ will denote the subset of $S_{\frac{3}{2}}(N, \chi)$ with the property that its elements map to cusp forms under the Shimura lifting. It is known that all the elements of $S_{\frac{3}{2}}(N, \chi)$, lying in the orthogonal complement of the space of theta series, lie in $S^*_{\frac{3}{2}}(N, \chi)$ (the orthogonal complement is with respect to the Peterson product).

### 2.3. Operators.

• Hecke operators for half integral weight forms: we shall only employ operators of the form $T(p^2)$, where $p$ runs through the set of natural primes. The operators of the form $T(p)$, $(p|N)$ as defined in (J.-P. Serre and H.M. Stark, 1976) will not be used.
• Other operators.

Let $f \in M_k(N, \chi)$ be a modular form of integral weight $k$, with Fourier coefficients $a(n)$. We shall use the following well known operators.

1. The Fricke involution $W_n$: $(f|W_n)(z) = N^{-k/2} z^{-k} f(\frac{-1}{Nz})$. We have $f|W_n \in M_k(N, \overline{\chi})$.

2. The shift operator $V_l$ ($l \in \mathbb{N}$): $(f|V_l)(z) = f(lz)$; $f|V_l \in M_k(lN, \chi)$.

2.4. BLOCKS OF COEFFICIENTS AND ASSOCIATED CHARACTERS.

Let $f$ be a modular form of half integral weight with Fourier coefficients $a(n)$ and $d$ a square-free positive integer. The block of coefficients corresponding to $d$ is the sequence of numbers $\{a(d), a(4d), a(9d), \ldots, a(n^2 d), \ldots\}$. If $k$ is an integer, $\psi$ a character (mod $N$), we define (see (G. Shimura, 1973), §3), the character $\psi_d$ (mod $dN$), by

$$\psi_d(m) = \psi(m) \left(\frac{-1}{m}\right)^k \left(\frac{d}{m}\right).$$

Our quadratic residue symbol is the same as the one defined in (G. Shimura, 1973), and also, we shall assume as in (G. Shimura, 1973) that $\chi(-1) = 1$ (for odd characters $\chi$, the space $S_{k+\frac{1}{2}}(N, \chi)$ consists of the zero form only).

## 3. Results

*Theorem A.* Let $k \geq 1$ be a natural integer and let $f \in S_{k+\frac{1}{2}}(N, \psi)$. When $k = 1$ we require that $f \in S_{\frac{3}{2}}^*(N, \chi)$. Consider a square free natural integer $d$. We shall make the following assumption:

(*) $f$ is an eigenform for all Hecke operators $T(p^2)$ for all prime factors $p$ of $N$ not dividing the conductor of $\psi_d$.

If the block of coefficients $\{a(d), a(4d), \ldots, a(n^2 d), \ldots\}$ is periodic in $n$, then it must vanish completely.

The proof consists of two steps. First, we prove the result under more restrictive conditions. This is accomplished in theorem 1. Second, we reduce the general case to the particular case discussed in theorem 1. We shall proceed with the statement and proof of theorem 1. By an irreducible natural number we mean 1 or a prime number.

*Theorem 1.* We are under the hypothesis in the statement of Theorem A. Moreover, assume that $\psi_d$ is a primitive character and there exists an irreducible number $l$ ($l \in \mathbb{N}$), such that

$$a(dn^2) = a(dm^2) \text{ for all } m \equiv n \pmod{l}.$$

Then $a(dn^2) = 0$ for all $n \geq 1$.

**P**roof:
The plan of the proof is as follows. By Shimura's lifting theorem there is a modular form of weight $2k$ whose associated $L$-series is the product of two Dirichlet series, one of them being a zeta function twisted by a character and the other one is formed with the coefficients $a(n^2)$ of $f$. We break the later one as a sum of zeta functions twisted by various characters modulo $l$. The terms corresponding to even characters multiply the common $L$-factor to give the $L$-series of some weight $k$-modular forms. We multiply everything by $\Gamma(s)$ and use the functional equation. An incompatibility between the appearing $\Gamma$-factors shows that it is impossible to have $L$-series with odd characters entering the initial equality and then we will be left with an equality of forms of different weights, which is also impossible unless all the coefficients are zero.

Now we proceed with the proof. We shall use Shimura's lifting theorem (G. Shimura, 1973), to find a modular form $F$ in $S_{2k}(N, \psi_1)$, whose Fourier expansion is

$$F(z) = \sum_{n \geq 1} A(n) q^n$$

where the coefficients $A(n)$ are defined via the series identity

$$\sum_{n \geq 1} A(n) n^{-s} = \left(\sum_{m \geq 1} \psi_d(m) m^{k-1-s}\right) \cdot \left(\sum_{m \geq 1} a(dm^2) m^{-s}\right)$$

and $\psi_d$ is the character defined in §2.4.
Using our hypothesis of periodicity for the $a(dn^2)$'s we can write the second factor appearing in the above series equality as

$$\sum_{m \geq 1} a(dm^2) m^{-s} =$$
$$= a(dl^2) \cdot \sum_{l | m} m^{-s} + \sum_{h \in (\mathbb{Z}/l\mathbb{Z})^\times} a(dh^2) \cdot \sum_{m \equiv h \pmod{l}} m^{-s} =$$
$$= a(dl^2) l^{-s} \zeta(s) + \sum_{h \in (\mathbb{Z}/l\mathbb{Z})^\times} a(dh^2) \cdot \zeta_{h,l}(s).$$

Here we used the notation $\zeta_{h,l}(s)$ for the classical partial zeta function

$$\zeta_{h,l}(s) = \sum_{m \equiv h \pmod{l}} m^{-s}.$$

Our goal is to use a convenient functional equation for $\sum_{m \geq 1} a(m^2) m^{-s}$ and in order to do so we shall rewrite each $\zeta_{h,l}$ as a linear combination of Dirichlet $L$-series associated to various Dirichlet characters modulo $l$.

It is known that

$$L(s, \chi) = \sum_{h=1}^{l-1} \chi(h) \cdot \zeta_{h,l}(s)$$

for all the characters $\chi$ of $(\mathbf{Z}/l\mathbf{Z})^\times$, (Washington, 1997).
This gives a set of $l-1$ equation in the unknowns $\zeta_{h,l}$ and its matrix is invertible hence we can write for each $h$

$$\zeta_{h,l}(s) = \sum_{\chi \in \widehat{(\mathbf{Z}/l\mathbf{Z})^\times}} \beta_\chi(h) \cdot L(\chi, s), \text{ where } \beta_\chi(h) \in \mathbf{C}.$$

One then obtains

$$\sum_{m \geq 1} a(dm^2) m^{-s} = a(dl^2) \cdot l^{-s} \cdot \zeta(s) + \sum_{\chi \in \widehat{(\mathbb{Z}/l\mathbb{Z})^\times}} \left( \sum_{h \in (\mathbf{Z}/l\mathbf{Z})^\times} a(dh^2) \beta_\chi(h) \right) \cdot L(s, \chi).$$

This is a good point to notice that for the trivial character $\mathbf{1} \in \widehat{(\mathbb{Z}/l\mathbb{Z})^\times}$, the $L$ function $L(s, \mathbf{1})$ is $(1 - l^{-s})\zeta(s)$, so if we put

$$\alpha_\chi = \sum_{h \in (\mathbb{Z}/l\mathbb{Z})^\times} a(h^2) \beta_\chi(h)$$

and collect the like terms, then the identity from Shimura's theorem becomes

$$\sum_{n \geq 1} A(n) n^{-s} = ((a(dl^2) - \alpha_{\mathbf{1}}) l^{-s} + \alpha_{\mathbf{1}}) \zeta(s) \cdot L(s - k + 1, \psi_d) +$$

$$+ \sum_{\chi \in \widehat{(\mathbf{Z}/l\mathbf{Z})^\times}, \chi \neq \mathbf{1}} \alpha_\chi \cdot L(s, \chi) \cdot L(s - k + 1, \psi_d)$$

The key of the proof is that half of these products of $L$-series give $L$ series corresponding to modular forms of weight $k$. More precisely, if $\chi$ is an even nontrivial character modulo $l$, then

$$(\psi_d \cdot \chi)(-1) = \psi_d(-1) = (-1)^k,$$

thus by theorem 4.7.1, pp 177 of [8] we know that the $L$-series product

$$L_{f^\chi}(s) = L(s,\chi) \cdot L(s-k+1, \psi_d)$$

gives the $L$-series associated to a modular form $f^\chi$ in $M_k(Nl, \psi_d\chi)$. A similar statement holds true for the trivial character $\mathbf{1}$. The $L$-series product

$$L_{f^{\mathbf{1}}}(s) = \zeta(s) \cdot L(s-k+1, \psi_1)$$

is the $L$-series associated to a modular form $f^{\mathbf{1}}$ in $M_k(N, \psi_1)$.

The hypothesis of the result quoted above requires $\psi_d$ to be a primitive character and this is why we had to impose this condition in the hypothesis of our theorem. The fact that the nontrivial characters $\chi$ of $(\mathbb{Z}/l\mathbb{Z})^\times$ are primitive $\pmod{l}$ is immediate from the fact that $l$ is a prime.

We collect all the terms corresponding to even characters and rewrite the above identity as

$$L_F(s) - ((a(l^2) - \alpha_{\mathbf{1}})l^{-s} + \alpha_{\mathbf{1}}) \cdot L_{f^{\mathbf{1}}}(s) - \sum_{\chi \text{ even}, \chi \neq \mathbf{1}} \alpha_\chi \cdot L_{f_\chi}(s) =$$
$$= \sum_{\chi \text{ odd}} \alpha_\chi L(s-k+1, \psi_d) \cdot L(s, \chi).$$

At this point we need to notice that $l^{-s} L_{f^{\mathbf{1}}}(s) = L_{f^{\mathbf{1}}|V(l)}(s)$ where $V(l)$ is the $l$-shift operator for modular forms and this way the left hand side of our equality is a linear combination of $L$-series associated to various modular forms.

Now we multiply both sides by $\Gamma(s)$ and notice that the right hand side is equal to

$$\Gamma\left(\frac{s}{2}\right) \cdot \Gamma\left(\frac{s-k+1+\delta}{2}\right)^{-1} \cdot \sum_{\chi \text{ odd}} \alpha_\chi \Lambda(s-k+1, \psi_d) \cdot \Lambda(s, \chi)$$

where, due to our condition on the parity of $\chi$ we put

$$\Lambda(s, \chi) = \Gamma((s+1)/2) L(s, \chi)$$

$$\Lambda(s, \psi_d) = \Gamma((s+\delta)/2) L(s, \psi_d),$$

with $\delta \in \{0, 1\}$ satisfying the congruence $k + \delta \equiv 0 \pmod{2}$.

The sum

$$\sum_{\chi \text{ odd}} \alpha_\chi \Lambda(s-k+1, \psi_d) \cdot \Lambda(s, \chi)$$

can be analytically continued to the whole plane while the product of the two gamma factors introduce infinitely many zeroes to the left on the real line. In fact all the negative odd integers sufficiently large are zeroes for this product.

This imposes strong conditions on the left hand side. In fact using the functional equation for each form entering the left hand side we will show that there exists a series (convergent for $\text{Re}(s) \gg 0$),

$$\sum_{n \geq 1} B(n) n^{-s}$$

with infinitely many zeroes to the right and the $B(n)$'s grow polynomially in $n$. But this is impossible unless all the $B(n)$'s are identically zero as it is proved in the appendix lemma. Then our initial series must be zero and this amounts to the assertion that a modular form of weight $2k$ is equal to a linear combination of modular forms of weight $k$. This forces both forms to be zero, hence all the $a(n)$'s must be zero.

Now we only need to show how the functional equations imply the desired conclusion. The term

$$\Lambda_F(s) = (\frac{2\pi}{\sqrt{(lN)}})^{-s} \Gamma(s) L_F(s)$$

has a functional equation of the form

$$\Lambda_F(s) = i^{2k} \Lambda_G(2k - s)$$

where

$$G(z) = F_{|\omega}, \ \omega = \begin{bmatrix} 0 & 1 \\ -lN & 0 \end{bmatrix}$$

and $G$ is a cusp form so its coefficients have a polynomial growth in $n$. On the other hand, for a modular form $f$ of weight $k$ whose $L$-series comes from products of Dirichlet $L$-series, one already has a functional equation of the form

$$\Lambda_f(s) = i^k \Lambda_g(k - s)$$

where

$$\Lambda_f(s) = (\frac{2\pi}{\sqrt{(lN)}})^{-s} \Gamma(s) L_f(s), \ \ g(z) = f_{|\omega}$$

and the coefficients have the desired growth as it is easy to see from formulas (4.7.10) and (4.7.15) in proof of Theorem 4.7.1 in [8]. Let $f$ denote the linear combination of all these weight $k$- forms entering the left hand side. If

$$\Lambda_F(s) - \Lambda_f(s)$$

has infinitely many zeroes to the left, then changing variable $s \longmapsto 2k - s$ one can see that

$$\Lambda_G(s) - \Lambda_g(s - k)$$

has infinitely many zeroes to the right. This is the desired series to which one can apply the lemma in order to force its vanishing. By symmetry the initial function must be zero and this forces

$$F = f$$

which is impossible unless both are zero.

<div style="text-align:right">Q.E.D.</div>

**P**roof of theorem A:

The general case where we drop the assumptions on the character $\psi_d$ and the period of the block $\{a(n^2 d) : n \in \mathbb{N}\}$ is similar to the one just proved above. By the same argument one can write the partial zeta function

$$\zeta_{h,M}(s) = \sum_{\chi \in \widehat{(\mathbf{Z}/M\mathbf{Z})^\times}} \alpha_\chi L(s, \chi) \quad (h \text{ and } M \text{ are coprime})$$

and the main problem we need to deal with is writing the $L$-series corresponding to non primitive characters in terms of ones which correspond to primitive characters, a problem which did not appear when the period $l$ was assumed to be a prime number or 1.
This problem is easy to overcome since any Dirichlet $L$-series can be reduced to one with a primitive character times some Euler factors. More precisely, let $\chi$ be a character of the multiplicative group $\widehat{(\mathbf{Z}/M\mathbf{Z})^\times}$ of conductor $\mathfrak{f}$, and let $\chi_0$ be its restriction to $\widehat{(\mathbf{Z}/\mathfrak{f}\mathbf{Z})^\times}$. Then it is known that

$$L(s, \chi) = \prod (1 - \frac{\chi_0(p)}{p^s}) \cdot L(s, \chi_0)$$

where the product is taken over all those primes $p$ dividing $M$ and prime to the conductor $\mathfrak{f}$. The argument goes through, as multiplying an $L$-series associated to a modular form by the extra factor $p^{-s}$ gives an $L$-series of the $p$-shift of the initial modular form.

<div style="text-align:right">Q.E.D.</div>

**R**emark: Our method can be used to prove statements of a similar nature, implying non-periodicity of sequences of the form $\{\frac{a(dn^2)}{n^i} | n \in$ $\mathbb{N}$}, where $d$ is a square-free positive integer and $i$ must satisfy the condition $i \neq k$. The only difference is that the argument $s$ of the $L$-series resulting from the break of $\sum_{n \geq 1} a(dn^2)n^{-s}$ is replaced by $s - i$. That changes the weight of the modular forms appearing as products of Dirichlet $L$-series by increasing it by $i$.

## 4. Appendix

In this appendix we shall give a proof for a lemma used in the above proof. It is known and it is part of the mathematical folklore but we add it here for completeness.

*Lemma 1.* Let
$$L(s) = \sum_{n \geq 1} B(n) n^{-s}$$
be a series with coefficients of polynomial growth in $n$,  $B(n) = O(n^\lambda)$. If as a function in $s$ this series has zeroes in each half plane $Re(s) > \sigma$, then it is zero.

**P**roof:

Let $B(n_0)$ be the first nonzero coefficient of this series. Then
$$L(s) = \frac{1}{n_0{}^s} \cdot \left( B(n_0) + \sum_{n \geq n_0} B(n) \cdot (n_0/n)^s \right).$$

Assume that $B(n) < Cn^\lambda$ for some positive constant $C$ and all $n \geq n_1$. Then we can split the sum between brackets into two sums and majorize its absolute value by
$$\sum_{n \geq n_0} |B(n)| \cdot (n_0/n)^{Re(s)} = \sum_{n_0 \leq n \leq n_1} |B(n)| \cdot (n_0/n)^{Re(s)} +$$
$$+ \sum_{n \geq n_1} n_0^\lambda \cdot (n_0/n)^{Re(s)-\lambda}.$$

It is now obvious that when $Re(s) >> 0$ the above two sums can be made as small as we want and this shows that when $Re(s)$ is large enough $L$ will not have any zero.

Q.E.D.

## Acknowledgements

The work for this paper was done at Johns Hopkins University and is part of author's Ph.D. thesis. The author wishes to express his thanks to Professor V. Kolyvagin for his suggestions and warm encouragement.